\documentclass[reqno,11pt]{amsart}

\usepackage{epsf}
\usepackage{graphics}
\usepackage{amssymb}
\usepackage{amsmath}
\usepackage{wasysym}

\date{}

\theoremstyle{plain}
\newtheorem{theorem}{Theorem}

\newtheorem{lemma}{Lemma}
\newtheorem{proposition}{Proposition}

\theoremstyle{definition}

\theoremstyle{remark}

\newtheorem*{examples}{Examples}
\newtheorem*{remark}{Remark}

\def\N{{\mathbb N}}
\def\Z{{\mathbb Z}}
\def\Q{{\mathbb Q}}
\def\R{{\mathbb R}}
\def\boA{\mathcal{A}}
\def\boB{\mathcal{B}}
\def\boC{\mathcal{C}}
\def\boD{\mathcal{D}}
\DeclareMathOperator{\wh}{Wh}
\DeclareMathOperator{\sub}{Sub}
\DeclareMathOperator{\aut}{Aut}
\DeclareMathOperator{\hair}{Hair}
\DeclareMathOperator{\loc}{Loc}
\DeclareMathOperator{\rat}{Rat}

\title{Asymptotic Vassiliev Invariants for Vector Fields}

\author{Sebastian Baader \and Julien March\'e} \thanks{first author supported by Fondation des Sciences Math\'ematiques de Paris}

\begin{document}

\maketitle

\begin{abstract} We analyse the asymptotical growth of Vassiliev invariants on non-periodic flow lines of ergodic vector fields on domains of $\R^3$. More precisely, we show that the asymptotics of Vassiliev invariants is completely determined by the helicity of the vector field. As an application, we determine the asymptotic Alexander and Jones polynomials and give a formula for the asymptotic Kontsevich integral.
\end{abstract}

\section{Introduction}

A smooth vector field on a manifold defines a flow whose orbits may be closed or not. If the manifold is a compact domain $G \subset \R^3$, we may ask about the asymptotical growth of knot invariants on non-periodic orbits. A well-known and classical example for this is the helicity of a vector field, which measures how pairs of non-periodic orbits are asymptotically linked, in the average~\cite{AK}. In order to make quantitative statements, we suppose that the flow of the vector field $X$ is measure-preserving and ergodic with respect to a probability measure $\mu$ on $G$, further, that the singularities of $X$ are isolated and that the periodic orbits of $X$ are not charged by the flow. For every non-periodic point $p \in G$ and $T>0$ we define a set $K(p,T) \subset \R^3$, as follows:
$$K(p,T)=\{\phi^X(p,t)| t \in [0,T]\} \cup [p, \phi^X(p,T)],$$
where $\phi^X$ is the flow of $X$ and $[p, \phi^X(p,T)]$ the geodesic segment in $\R^3$ joining $p$ and $\phi^X(p,T)$. This set is actually a knot, i.e. an embedded circle, for almost all $p \in G$, $T>0$ (\cite{GG1}, \cite{V}). Under the above hypotheses, Gambaudo and Ghys proved the existence of an asymptotic signature invariant which is proportional to the helicity of $X$~\cite{GG1}: for almost all $p \in G$ the limit
$$\sigma(X)=\lim_{T \rightarrow \infty} \frac{1}{T^2} \sigma(K(p,T)) \in \R$$
exists and is independent of the starting point $p \in G$. Here $\sigma$ denotes the signature invariant of links. The asymptotic signature invariant determines the asymptotical behaviour of a large class of concordance invariants~\cite{Ba}. In this note we show that Vassiliev invariants are asymptotically determined by the signature (hence also by the helicity).

\begin{theorem} \label{main} Let $v$ be a real-valued Vassiliev knot invariant of degree $n$. There exists a constant $\alpha_v \in \R$, such that for almost all $p \in G$ 
$$\lim_{T \rightarrow \infty} \frac{1}{T^{2n}} v(K(p,T))=\alpha_v \sigma(X)^n.$$
The constant $\alpha_v$ does not depend on the vector field $X$.
\end{theorem}

Gambaudo and Ghys provided the first instance of this theorem since the helicity can be defined as an asymptotical linking number, which is a Vassiliev invariant of degree one (for links, however). The proof of Theorem~\ref{main} is based on a asymptotical count of Gauss diagrams with respect to suitable diagrams of the knots $K(p,T)$. We give a short summary of Gambaudo and Ghys' construction in Section~2. The proof of Theorem~\ref{main} is contained in Section~3. The rest of the paper is devoted to applications. In Section~4, we determine the asymptotic Alexander and Jones polynomials and give a formula for the asymptotic Kontsevich integral, using the technique of wheels and wheelings. We conclude with an appendix where we give an expression for the Kontsevich integral of torus knots, simplifying the arguments of \cite{Ma}. From this, we deduce an explicit tree expansion for the asymptotic Kontsevich integral.

\begin{remark} For reasons of simplicity, we restrict ourselves to the study of asymptotical knots and their invariants, rather than links. The case of links does not pose any additional difficulties.
\end{remark}

\section{Asymptotic Diagrams}

The main part of Gambaudo and Ghys' work~\cite{GG1} consists in constructing good diagrams for the knots $K(p,T)$. For this purpose, they cover the domain $G$, away from the singularities of $X$, by a countable family of flow boxes $\{\mathcal{F}_i\}_{i \in \N}$. Further, they define a projection $\pi: \R^3 \rightarrow \R^2$ onto a plane which is well-adapted to this family: for every $\epsilon>0$, there exists a finite subset $\mathcal{C} \subset \N$, such that for almost all $p \in G$, $T>0$ large enough, the diagram $\pi(K(p,T))$ is regular and, up to an error $\leq \epsilon T^2$, its crossings arise from pairs of overcrossing flow boxes $\mathcal{F}_i$, $\mathcal{F}_j$, with $i,j \in \mathcal{C}$. Moreover, at these finitely many overcrossing spots the diagram looks like a rectangular grid, as sketched in Figure~1. We will shortly see that the number of crossings of these grids grows like $T^2$.

\begin{figure}[ht]
\scalebox{1}{\raisebox{-0pt}{$\vcenter{\hbox{\epsffile{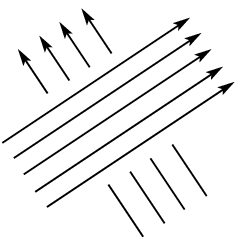}}}$}} 
\caption{}
\end{figure}

Let $n_i(p,T)=\pi_0(\mathcal{F}_i \cap \{\phi^X(p,t)| t \in [0,T]\})$ be the number of times the flow line starting at $p$ and ending at $\phi^X(p,T)$ enters the flow box $\mathcal{F}_i$. Applying Birkhoff's ergodic theorem to the characteristic function of the flow box $\mathcal{F}_i$, we immediately see that for almost all $p \in G$ the limit
$$n_i=\lim_{T \rightarrow \infty} \frac{1}{T} n_i(p,T)>0$$
exists (and is proportional to the volume of the flow box $\mu(\mathcal{F}_i)$). Therefore the number of crossings $c_{ij}(p,T)$ at an overcrossing spot of two flow boxes $\mathcal{F}_i$, $\mathcal{F}_j$ satisfies
\begin{equation}
\lim_{T \rightarrow \infty} \frac{1}{T^2} c_{ij}(p,T)=n_i n_j.
\label{average1}
\end{equation}
For later purposes, we choose a natural number $N \in \N$ and subdivide the time interval $[0,T]$ into $N$ sub-intervals $I_1, I_2, \ldots, I_N$ of length $\frac{T}{N}$. Every index $k \in \{1,2,\ldots, N \}$ gives rise to a function
$n_{i,k}(p,T)=\pi_0(\mathcal{F}_i \cap \{\phi^X(p,t)| t \in I_k\})$. Again, by Birkhoff's theorem, we obtain
\begin{equation}
\lim_{T \rightarrow \infty} \frac{1}{T} n_{i,k}(p,T)=\frac{n_i}{N}.
\label{average2}
\end{equation}
At last, for two flow box indices $i_1,i_2 \in \mathcal{C}$, and $k_1, k_2 \in \{1,2,\ldots, N \}$, we define the number of crossings $c_{i_1,k_1,i_2,k_2}(p,T)$ between $\mathcal{F}_{i_1} \cap \{\phi^X(p,t)| t \in I_{k_1}\})$ and $\mathcal{F}_{i_2} \cap \{\phi^X(p,t)| t \in I_{k_2}\})$ at an overcrossing spot of the flow boxes $\mathcal{F}_{i_1}$, $\mathcal{F}_{i_2}$. Equation~(\ref{average2}) implies
\begin{equation}
\lim_{T \rightarrow \infty} \frac{1}{T^2} c_{i_1,k_1,i_2,k_2}(p,T)=\frac{n_{i_1}n_{i_2}}{N^2}
\label{average3}
\end{equation}
This equality will play an important role in the proof of Theorem~\ref{main}.

\section{Gauss Diagram Formulae and Proof of Theorem~\ref{main}}

A Gauss diagram is nothing but a special notation for a knot diagram. It consists of an oriented circle with a finite number of signed arrows connecting pairs of points on the circle. The circle stands for the oriented knot itself, while the arrows encode crossing points of the knot diagram, pointing from the lower to the upper strand. Their signs indicate the signs of their crossings. For example, Figure~2 shows a Gauss diagram representing the standard diagram of the twist knot with six crossings. Here the orientation of the circle is understood to be clockwise.

\begin{figure}[ht]
\scalebox{1}{\raisebox{-0pt}{$\vcenter{\hbox{\epsffile{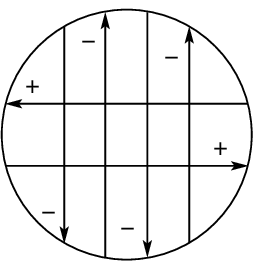}}}$}} 
\caption{}
\end{figure}

It is often convenient to consider pointed Gauss diagram, i.e. Gauss diagrams with a distinguished base point on its circle. Throughout this section, we will work with pointed Gauss diagram. In particular, we will be concerned with the pointed Gauss diagrams $G(p,T)$ arising from Gambaudo and Ghys' special diagrams $D(p,T)=\pi(K(p,T))$.

Gauss diagrams are of special interest in the theory of Vassiliev invariants, since the latter can be identified with certain formal linear combinations of Gauss diagrams. In order to explain this, we have to introduce a pairing between Gaus diagrams. Let $\Gamma$, $G$ be two Gauss diagrams. The expression
$$\langle \Gamma, G \rangle$$
is defined as the weighted number of sub-diagrams of $G$ isomorphic to $\Gamma$, respecting the circles, base points and all orientations. The weights are simply the products over all signs of arrows of the corresponding subgraphs and $\Gamma$. The above pairing evidently extends to a bilinear form on formal linear combinations of Gauss diagrams with real coefficients. A well-known theorem by Goussarov, Polyak and Viro~\cite{GPV} says that every real-valued Vassiliev invariant $v$ of order $n$ can be represented by a finite linear combination of Gauss diagrams $a_1 \Gamma_1+ \dots + a_m \Gamma_m$ with at most $n$ arrows: for every knot $K$, for every Gauss diagram $G$ representing $K$, the value $v(K)$ coincides with $\langle a_1 \Gamma_1+ \dots + a_m \Gamma_m, G \rangle$. However, not every linear combination of Gauss diagram gives rise to a knot invariant. 

\begin{examples} \quad
\begin{enumerate}
\item The pairing with the sum of the two pointed Gauss diagrams with one positive arrow does not define a knot invariant since it computes the writhe of a diagram. \\
\item The Casson invariant can be defined as the pairing with the diagram
$\scalebox{0.5}{\raisebox{-0pt}{$\vcenter{\hbox{\epsffile{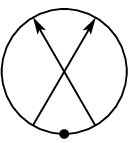}}}$}}$, as shown in~\cite{PV}. 
\\
\end{enumerate}
\end{examples}
The representation of Vassiliev invariants by linear combinations of Gauss diagrams allows us to reduce the proof of Theorem~\ref{main} to the following combinatorial statement.  

\begin{lemma} \label{pairing} Let $\Gamma$ be a Gauss diagram with $n$ arrows. For almost all $p \in G$ 
$$\lim_{T \rightarrow \infty} \frac{1}{T^{2n}} \langle \Gamma, G(p,T) \rangle=\frac{1}{(2n-1)!!}\sigma(X)^n.$$
\end{lemma}

\begin{remark} It is tempting to replace $\langle \Gamma, G(p,T) \rangle$ by $\langle \Gamma, K(p,T) \rangle$ in the above statement. This makes no sense unless the pairing with $\Gamma$ defines a knot invariant. Nevertheless, Theorem~\ref{main} is an immediate consequence of Lemma~\ref{pairing}, by the theorem of Goussarov, Polyak and Viro. 
\end{remark}

\begin{proof}[Proof of Lemma~\ref{pairing}]
Every occurence of $\Gamma$ in the diagram $G(p,T)$ corresponds to a collection of $n$ crossings involving $2n$ strands whose order is prescribed by $\Gamma$. We start by observing that only the crossings arising from overcrossing flow boxes $\mathcal{F}_i$, $\mathcal{F}_j$, with $i,j \in \mathcal{C}$, produce an essential contribution to $\langle \Gamma, G(p,T) \rangle$. Indeed, the total number of $n$-tuples of crossings of $D(p,T)$ grows like $T^{2n}$, by equation~(\ref{average1}). If we restrict ourselves to $n$-tuples with at least one `exceptional' crossing, then this number is estimated by a constant multiple of $\epsilon T^{2n}$. 

Let us now choose a natural number $N \in \N$ and subdivide the interval $[0,T]$ into $N$ sub-intervals $I_1, \ldots, I_N$, as in Section~2. Every $n$-tuple of crossings $(c_1, \ldots, c_n)$ of $D(p,T)$ defines $2n$ times $t_1, \ldots, t_{2n} \in [0,T]$, with the convention that $t_{2k-1}$, resp. $t_{2k}$, denote the first, resp. second, time of occurence at the crossing $c_k$ ($k \in \{1, \ldots, n \}$), starting at $t=0$. The order of the $t_l$'s tells us whether the $n$-tuple contributes to $\langle \Gamma, G(p,T) \rangle$ or not. We would like to reduce our count to $n$-tuples with sufficiently distinct times, i.e. such that every interval $I_j$ contains at most one $t_k$ (this naturally forces $N \geq 2n$). Let us fix two flow box indices $i_1,i_2 \in \mathcal{C}$ and determine $c_{i_1,j,i_2,j}(p,T)$, the number of (single) crossings between $\mathcal{F}_{i_1}$, $\mathcal{F}_{i_2}$ for which both times belong to the same interval $I_j$, $j \in \{1, \ldots N\}$. By equation~(\ref{average3}), this number grows like $\frac{n_{i_1}n_{i_2}}{N^2}T^2$. Since there are $N$ intervals $I_j$, we have to sum up the $c_{i_1,j,i_2,j}(p,T)$ over $j$ and conclude that the number of single crossings between $\mathcal{F}_{i_1}$, $\mathcal{F}_{i_2}$ for which both times belong to the same interval grows like $\frac{n_{i_1}n_{i_2}}{N}T^2$. Choosing $N$ large enough, we see that these crossings do not produce an essential contribution to $\langle \Gamma, G(p,T) \rangle$, either.

It remains to count $n$-tuples $c=(c_1, \ldots, c_n)$ of crossings where all $t_l$'s belong to different intervals $I_j$, whose order is prescribed by $\Gamma$. These $n$-tuples come together with injective maps $f_c: \{1, 2, \ldots 2n\} \rightarrow \{1,2, \ldots N\}$ satisfying $t_l \in I_{f_c(l)}$. Let us fix such a map $f$ which is compatible with the Gauss diagram $\Gamma$, meaning that the obvious Gauss diagram associated with $f$ coincides with $\Gamma$. Moreover, let us fix $n$ not necessarily distinct spots $X_1, \ldots, X_n$ of overcrossing flow boxes. We denote by $i_k, j_k \in \mathcal{C}$ the indices of the two flow boxes that cross at $X_k$. Using equation~(\ref{average3}) again, we see that the number of $n$-tuples $c$ with prescribed function $f_c$ and $c_k \in X_k$ grows like
$$\prod_{k=1}^n \frac{n_{i_k}n_{j_k}}{N^2}T^2.$$ 
Summing over all injective maps $f$ compatible with $\Gamma$ and over all collections of $n$ spots $X_k$ (finitely many in number), we obtain an expression that grows like $T^{2n}$. However, this is not quite $\langle \Gamma, G(p,T) \rangle$, since we have been neglecting signs so far. This will finally lead us to the desired relation with the asymptotic signature.

We observe that all crossings of the grid at a spot $X_k$ share the same sign $\epsilon_k$. Taking these into account, we obtain the following expression for the asymptotical growth of $\langle \Gamma, G(p,T) \rangle$:
\begin{equation}
\sum_f \sum_{\{X_k\}} \prod_{k=1}^n \epsilon_k n_{i_k} n_{j_k} \frac{T^2}{N^2},
\label{sum}
\end{equation}
where the first sum is taken over all injective maps $f$ compatible with $\Gamma$ and the second one over all collections of $n$ spots $X_k$. The inner sum is easily recognizable as
$$\left(\sum_{X_k} \epsilon_k n_{i_k} n_{j_k} \frac{T^2}{N^2} \right)^n,$$
where the sum is now taken over all single spots $X_k$. For large $N \in \N$, the proportion of injective maps $f: \{1, 2, \ldots 2n\} \rightarrow \{1,2, \ldots N\}$ compatible with $\Gamma$ converges to $\frac{n!}{(2n)!}$, i.e. their number grows like $\frac{n!}{(2n)!} N^{2n}$. To see this, we note that the asymptotic proportion will not depend on the diagram $\Gamma$, since we can reorder the intervals $I_1, \ldots, I_N$. Thus we need only consider the case of the Gauss diagram corresponding to strictly monotone maps $f$. The proportion of these converges to
$$\int_{0 \leq t_1 \leq \ldots \leq t_{2n} } dt_1 \ldots dt_{2n}=\frac{1}{(2n)!}.$$
Then, there are $n!$ permutation of pairs of times $(t_{2k-1},t_{2k})$ which do not affect the Gauss diagram. This leads to the asymptotic proportion $\frac{n!}{(2n)!}$.
The expression~(\ref{sum}) therefore simplifies to
$$\frac{n!}{(2n)!} \left(\sum_{X_k} \epsilon_k n_{i_k} n_{j_k} T^2 \right)^n.$$
The sum between brackets is nothing but the asymptotical writhe of the diagram $D(p,T)$, which converges to twice the asymptotical signature of $K(p,T)$, by Gambaudo and Ghys~\cite{GG1}. At last, we obtain
$$\lim_{T \rightarrow \infty} \frac{1}{T^{2n}} \langle \Gamma, G(p,T) \rangle=\frac{n!}{(2n)!}(2 \sigma(X))^n=\frac{1}{(2n-1)!!}\sigma(X)^n.$$
\end{proof}

\section{Asymptotic Kontsevich Integral}

\subsection{Reduction to Torus Knots}

Let $X$ be a vector field on a bounded domain of $\R^3$ as in the settings of Theorem~\ref{main} and denote by $\sigma$ its asymptotic signature invariant. Then, for any finite type invariant $v$ of degree $n$, the constant $\alpha_v$ describes the asymptotics of $v(K(p,T))$ for $T$ going to infinity. Precisely, we denote by $\alpha_n(v)$ the limit of $\frac{v(K(p,T))}{T^{2n}\sigma^n}$, where $\sigma$ is supposed to be non-zero.

The correspondence $v\mapsto \alpha_n(v)$ is linear and if $v$ happens to have degree less than $n$, then $\alpha_n(v)=0$. We conclude that $\alpha_n$ is a linear form on the $n$-th graded space of Vassiliev invariants which itself is dual to the linear space of $n$-chord diagrams $\boA_n(S^1)$. Hence, $\alpha_n$ belongs to $\boA_n(S^1)$ and our purpose in this section is to compute $\alpha=\sum_n \alpha_n\in \boA(S^1)$.

The trick is to reduce the computation to torus knots, thanks to the following lemma. Denote by $T(p,q)$ the torus knot with parameters~$p$ and~$q$.

\begin{lemma} \label{torus}
Let $\lambda$ be an irrational number in $[0,1]$. 
Then, there exist two sequences $p_k,q_k$ of coprime integers going to infinity such that $q_k/p_k$ converges to $\lambda$  and $v(T(p_k,q_k))/{p_k}^{2n}$ converges to  $\alpha_n(v)\frac{\lambda^n}{2^n}$, for every Vassiliev invariant $v$ of order $n$.
\end{lemma}

\begin{proof}
Consider the square $[0,1]\times[0,1]$ and identify two sides according to the rule $(x,0)=(x+\lambda \mod 1,1)$. It is well-known that the vector field $(0,1)$ defines an ergodic flow on this space. 
Cut the square along the segment $\{1-\lambda\}\times[1/2,1]$ and embed this template as in Figure~\ref{template}. We note that this template can be drawn on a standard torus in $\R^3$. The proof of Theorem~\ref{main} implies that for almost all starting points $p$ on the template, $v(K(p,T))/T^{2n}$ converges to $\alpha_n(v) \frac{\lambda^n}{2^n}$. 
Let $p=(x_0,y_0)$ be such a starting point: for all integers $k$, $\phi(p,k)=(x_k,y_0)$, where $x_k=x_0+k\lambda\mod 1$. We define a sequence of times $t_k$ for which the flow goes very close to the starting point. In formulas, set $t_0=0$ and $t_{k+1}=\inf \{m > t_k, x_0<x_m<\inf\{x_l, 0<l<m\}\}$. Denote by $K_k$ the knot obtained by closing the orbit of $p$ during a time $t_k$ with the segment $[x_0,x_{t_k}]\times\{y_0\}$. 

\begin{figure}[ht]
\scalebox{0.7}{\raisebox{-0pt}{$\vcenter{\hbox{\epsffile{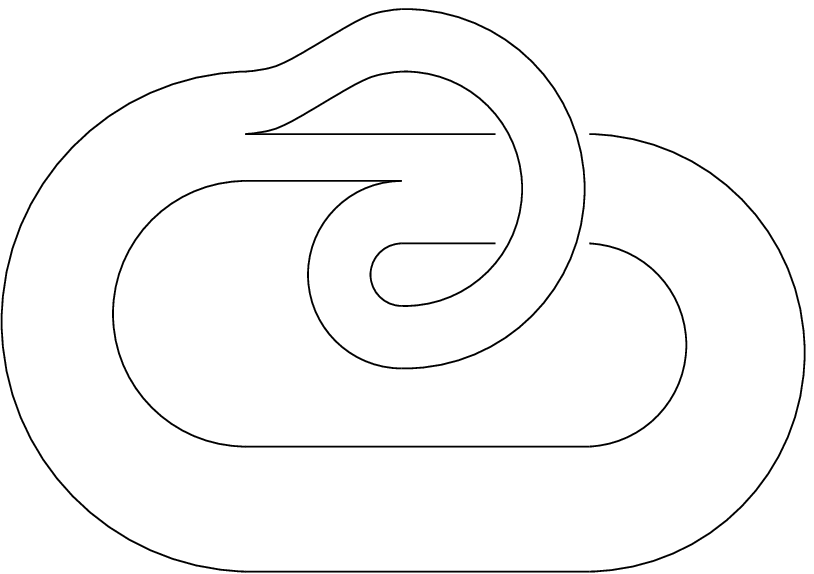}}}$}} 
\caption{} \label{template}
\end{figure}

This knot can be drawn on a torus and hits the meridian $t_k$ times. The number of times it hits the parallel is 
$[\lambda t_k]$. The two sequences $p_k=t_k$ and $q_k$ satisfy the assumptions of the lemma.
\end{proof}

As an easy application, one can give formulas for the asymptotic Alexander polynomial and Jones polynomial.

\begin{proposition}
Let $X$ be a vector field as in Theorem~\ref{main} with asymptotic signature invariant $\sigma$.
 
Set $\tilde{\Delta}=\Delta(e^h)\in \Q[[h]]$, where $\Delta$ is the Alexander polynomial normalized such that it is a symmetric Laurent polynomial.
$$\lim_{T\to\infty}\tilde{\Delta}_{K(p,T)}(h/T^2)=\frac{\sinh(\sigma h)}{\sigma h}$$

Set $\tilde{J}=J(e^h)\in \Q[[h]]$, where $J$ is the Jones polynomial normalized such that it takes the value 1 at the trivial knot.
$$\lim_{T\to\infty}\tilde{J}_{K(p,T)}(h/T^2)=-h\sigma e^{\sigma h}$$
\end{proposition}

\begin{proof}
This is an easy application of Lemma~\ref{torus} using the formulas 
$\Delta_{T(p,q)}(t)=t^{-(p-1)(q-1)/2}\frac{(t^{pq}-1)(t-1)}{(t^p-1)(t^q-1)}$ and $J_{K(p,q)}(t)=\frac{t^{(p-1)(q-1)/2}(1-t^{p+1}-t^{q+1}+t^{p+q})}{1-t^2}$. 
\end{proof}

\subsection{The Kontsevich Integral of Torus Knots and its Limit}

Let $\boA$ be the space of formal series of trivalent diagrams lying on $S^1$ up to the usual (AS) and (IHX) relations. The Kontsevich integral $Z(K)$ of a knot $K$ takes values in $\boA$. It happens to be easier to work with a space of formal series of uni-trivalent diagrams modulo the same relations denoted by $\boB$. We define the degree of a diagram as half the total number of vertices and the product of two diagrams as their disjoint union.

The wheeling map $\Upsilon:\boB\to\boA$ is a well-known isomorphism of algebras between these two spaces. We will denote by $Z^{\sun}\!(K)$ the series $\Upsilon^{-1}Z(K)\in \boB$: we will need a closed expression for $Z^{\sun}\!(T(p,q))$ for relatively prime integers $p$ and $q$.

Let us start with equation (2) of \cite{Ma}:
\begin{equation}\label{formule}
Z^{\sun}T(p,q)=\partial_{\Omega}^{-1}\left(_q\Omega\cdot\Omega_p\exp(\frac{pq}{2}\!\!\frown)\right)
\frac{\exp(-\frac{pq}{2}\!\!\frown+\frac{pq}{48}\Theta)}{\langle \Omega,\Omega\rangle}.
\end{equation}
In this formula, $\frown$ is the diagram of $\boB$ consisting of a single edge, $\Theta$ is the theta graph and $\Omega$ is a series of diagrams to be defined later. Given two diagrams $A$ and $B$ one defines the following operations:
\begin{enumerate}
\item $\langle A,B\rangle$ is the sum of all diagrams obtained by gluing all univalent vertices of $A$ with all univalent vertices of $B$.
\item $\partial_A B$ is the sum of all diagrams obtained by gluing all univalent vertices of $A$ with some univalent vertices of $B$.
\item $A\cdot B$ is the sum of all diagrams obtained by gluing some univalent vertices of $A$ with some univalent vertices of $B$.
\item $_qA\cdot B_p$ is the same sum as before except that each term is multiplied by $q^ap^b$, where $a$ and $b$ are the numbers of remaining univalent vertices of $A$ and $B$ respectively.
\end{enumerate}

Given a series $f(h)\in \R[[h]]$ we define the wheel series $\wh(f(h))\in \boB$ by associating to each monomial $h^n$ a circle with $n$ legs attached to it: for instance one has $\wh(h^8)=\sun$.
These series of diagrams play a central role: if we define $\Omega=\exp(\wh(F(h)))$, where $F(h)=\frac{1}{2}\log\frac{\sinh(h/2)}{h/2}$, we have for instance $Z^{\sun}(U)=\Omega/\langle\Omega,\Omega\rangle$, where $U$ is the unknot.

\begin{proposition} \label{limit}
Let $X$ be a vector field as in Theorem~\ref{main} with asymptotic signature invariant $\sigma$.
Given $\alpha\in\R$ and $D\in \boB$ we define $\alpha^{\deg}D$ as $\alpha^{n}D$, where $n$ is the degree of $D$. Then, the following formula holds:
$$\lim_{T\to\infty} T^{-2\deg}Z^{\!\sun} K(p,T)=(2\sigma)^{\deg}\partial_{\Omega}^{-1}\left( \exp(\frac{1}{2}\!\frown)\right) \exp(-\frac{1}{2}\!\frown+\frac{\Theta}{48}).$$
\end{proposition}

\begin{proof}
Thanks to Lemma~\ref{torus}, it is sufficient to consider the limit of $Z^{\!\sun} T(p,q)$ for $p$ and $q$ going to infinity at the same rate. The corresponding asymptotic signature will be equal to $1/2$. Thanks to Theorem~\ref{main}, we recover the general formula by applying the operator $(2\sigma)^{\deg}$.

Consider formally the expression $p^{-2\deg}Z^{\!\sun}T(p,p)$ in equation \eqref{formule}. 
The expression $p^{-2\deg}\langle\Omega,\Omega\rangle$ goes to 1 as any non empty diagram in $\langle\Omega,\Omega\rangle$ is divided by a positive power of $p$.

As $p^{-2\deg}$ commutes with addition and multiplication in $\boB$, one has $p^{-2\deg}\exp(-\frac{p^2}{2}\!\frown+\frac{p^2}{2}\Theta)=\exp(-\frac{1}{2}\!\frown+\frac{\Theta}{48})$.

Remark that $\partial_{\Omega}^{-1}=\partial_{\Omega^{-1}}$ and $\Omega^{-1}=\exp(-\wh(F(h)))$. Consider now a term of the expression $p^{-2\deg}\partial_{\Omega^{-1}}\left(_p\Omega\cdot\Omega_p\exp(\frac{p^2}{2}\!\!\frown)\right)$. It is a diagram $D$ obtained by gluing wheels coming from $\Omega^{-1}\!,\,_p\Omega$ and $\Omega_p$. Let us denote by $z,x,y$ these three distinct types of wheels appearing in $D$. The diagram $D$ is multiplied by a factor $p^{2n_z+n_x+n_y}$, where $n_z$ is the number of edges which are attached only to wheels of type $z$. The numbers $n_x$ (resp. $n_y$) denote the number of edges attached to wheels of type $x$ (resp. of type $y$) which are not attached to wheels of type $y$ (resp. of type $x$). 

We remark that the degree of $D$ is equal to the total number of edges attached to the wheels. Hence, $p^{-2\deg}D$ has a negative power of $p$ unless there are no wheels of type $x$ and $y$. When $p$ goes to infinity, only gluings between $\Omega^{-1}$ and $\exp(\frac{p^2}{2}\!\frown)$ remain and one obtains the desired formula.
\end{proof}

\appendix

\section{The Tree Expansion of the Kontsevich Integral}

\subsection{Gluing Graphs and their Substitutions}

Let $X$ be a finite set. We denote by $S(X)$ the vector space generated by isomorphism classes of finite graphs with vertices labeled by $X$ and cyclic orientation of the edges around vertices. We complete $S(X)$ with respect to the degree which counts edges and vertices.

Given $\Gamma\in S(X)$ and $f(\cdot,h):X\to \R[[h]]$, we define $\sub(\Gamma,f)$ in the following way: replace a vertex of $\Gamma$ labeled with $x$ by the wheel series $\wh(f(x,h))$ (using multilinearity), then sum over all possible gluings (respecting cyclic orientations) of the edges of $\Gamma$ to univalent vertices of wheels  located at the corresponding vertices. As a simple example, one has
$$\sub(\exp(\bullet),F(h))=\Omega.$$
We can interpret formula \eqref{formule} in terms of substitution of gluing graphs.
Set $X=\{x,y,z\}$ and denote by $G$ the following element of $S(X)$.
$$G=\sum_{[\Gamma]} \frac{(pq)^{-|E(\Gamma)|}}{|\aut (\Gamma)|} \,\Gamma.$$
The notation $[\Gamma]$ means that we sum over isomorphism classes of graphs $\Gamma$ which do not have edges between two vertices labeled by $x$ or two vertices labeled by $y$.

Set $f(x,h)=F(ph),f(y,h)=F(qh)$ and $F(z)=-F(pqh)$. Then, the definition of $G$ is designed to verify the following formula, taking care of the factors of $p$ and $q$ as in the proof of Proposition \ref{limit}.
\begin{equation}\label{sub}
\partial_{\Omega}^{-1}\left(_q\Omega\cdot\Omega_p\exp(\frac{pq}{2}\!\!\frown)\right)=\sub(G,f)\exp(\frac{pq}{2}\!\!\frown).
\end{equation}

Let us recall the notion of $\Pi$-diagram from \cite{Ma}. Let $\boC$ be the category whose objects are  free abelian groups of finite rank and morphisms are linear isomorphisms. If $\Pi$ is a functor from $\boC$ to real vector spaces, we define $\boD(\Pi)$ as a quotient of $\bigoplus_{[\Gamma]} \Pi(H^1(\Gamma,\Z))$ by generalized (AS) and (IHX) relations. In this direct sum, $\Gamma$ runs over trivalent diagrams.

As an example, if $\Pi(H)=\prod_nS^n(H\otimes\R)$, then there is an isomorphism between $\boB$ and $\boD(\Pi)$. Through this isomorphism, a leg is replaced by the cohomology class of the edge to which it is attached.
Setting $\Pi_s(H)=\Pi(H)[(H\otimes\R\setminus\{0\})^{-1}]$, one defines a space of singular diagrams where non-trivial legs become invertible. We denote this space by $\boB_s$. 

The interest of this space of diagrams comes from  Proposition 3.1 of \cite{Ma} that we recall without proof, though in a more explicit form.

\begin{proposition}\label{subs}
Let $\Gamma$ be a graph in $S(X)$ and $f(\cdot,h):X\to \R[[h]]$ a decoration of $X$ by formal series.
Denote by $\Gamma_{\!\circ}$ the trivalent diagram obtained from $\Gamma$ by gluing all edges incoming to a same vertex to a circle in the order given by the cyclic orientations on $\Gamma$.

For each vertex $v$ of $\Gamma$, we denote by $k_v$ the valency of $v$. Consider the edges of the circle lying at $v$ as cohomology classes in $H^1(\Gamma_{\!\circ},\Z)$. We introduce variables $y^1_v,\ldots, y^{p_v}_v$ to parametrize the $p_v$ distinct classes occuring with multiplicities $m^1_v,\ldots,m^{p_v}_v$. Then, the following formula holds:

$$\sub(\Gamma,f)=\prod_v \sum_{l_v=1}^{p_v}
\frac{\partial^{k_v-p_v}}{\partial^{m^1_v-1}_{y^1_v}\cdots\partial^{m^{p_v}_v-1}_{y^{p_v}_v}}
\left(
\frac
{
f'(v,y^{l_v}_v)
\prod_{j_v=1}^{p_v}
\frac{(y^{j_v}_v)^{m^{j_v}_v-1}}
{(m^{j_v}_v-1)!}
}
{
(y^{l_v}_v)^{k_v-p_v}
\prod_{j_v\ne l_v}(y^{l_v}_v-y^{j_v}_v)
}
\right)
$$
\end{proposition}

Let us summarize what is important in this formula:
\begin{itemize}
\item[-] there is a finite number of terms,
\item[-] the formula depends on $f'$ rather than $f$,
\item[-] all terms have non trivial denominator except if $p_v=1$ for all $v$. This occurs if and only if $\Gamma$ is a tree.
\end{itemize}

\subsection{Rationality}

Let $\Pi_{\rat}(H)=\R[\exp(H)]_{\loc}$, where we localize expressions of the form $P_1(e^{h_1})\cdots P_k(e^{h_k})$ for non-zero polynomials  $P_1,\ldots P_k$ and non-zero cohomology classes $h_1,\ldots,h_k$. Call $\boB_{\rat}=\boD(\Pi_{\rat})$.
The Taylor expansion $\Pi_{\rat}(H)\to\Pi(H)$ induces a map $\hair:\boB_{\rat}\to\boB_s$. Garoufalidis and Kricker (see \cite{GK}) showed that for every knot $K$, there exists a series $Z^{\sun}_{\rat}(K)\in \boB_{\rat}$ such that 
$$Z^{\sun}(K)=\frac{1}{\langle \Omega,\Omega\rangle}\exp\left(\wh(F(h)-\frac{1}{2}\log\tilde{\Delta}(h))\right)\hair Z^{\sun}_{\rat}(K).$$
Let us ignore from now on all wheel diagrams as they are well understood and are not in the image of the $\hair$ map. We will denote this part by $Z^{\sun}_{>1}$, as all non-wheel diagram have first Betti number greater than 1.

Recall that one has $Z^{\sun}T(p,q)=\frac{\exp(\frac{pq\Theta}{48})}{\langle \Omega,\Omega\rangle}\sub(G,f)$. The derivative of $F(h)$ is $\frac{1}{4}\frac{e^h+1}{e^h-1}-\frac{1}{2h}$, which is a mixture of a rational function in $e^h$ and a Laurent polynomial in $h$. In the substitution  of $G$ with $f$ there will appear a more complicated mixture which, at the end, happens to be rational, i.e. belongs to the image of the $\hair$ map. The point is that we can identify from the beginning which terms can be rational.

Consider $\Pi_{m}(H)$, the subalgebra of $\Pi_s(H)$ generated by $\Pi_{\rat}(H)$ and Laurent polynomials in $H$. This subalgebra is stable by derivation and Proposition \ref{subs} shows that the Kontsevich integral $Z^{\sun}_{>1}T(p,q)$ belongs to $\boD(\Pi_m)$ (except for the wheel part). The degree of a Laurent polynomial in $H$ extends to a well-defined degree on $\Pi_{m}(H)$ and on $\boD(\Pi_m)$ by functoriality. The rationality of the Kontsevich integral of torus knots shows that  $Z^{\sun}_{>1}T(p,q)$ belongs to the degree $0$ part.

Analysing Proposition~\ref{subs}, one sees that the degree 0 part is obtained as follows. 
Set $F^{+}(h)=\frac{1}{2}\log \sinh(h/2)$. It is not a Laurent power series but its derivative is. Setting $f^+(x,h)=F^+(ph),f^+(y,h)=F^+(qh),f^+(z)=-F^+(pqh)$, the formula of Proposition \ref{subs} allows us to make sense of $\sub(G,f^+)$. Then $\sub(G^t,f^+)$ is the degree 0 contribution to $Z^{\sun}T(p,q)$, where $G^t$ is the tree part of $G$. We deduce the following formula.

\begin{proposition}[The rational expression of the Kontsevich integral of torus knots]
\begin{equation}\label{ntor}
Z^{\sun}_{>1}T(p,q)=\frac{\exp(\frac{pq\Theta}{48})}{\langle\Omega,\Omega\rangle}\sum_{[\Gamma,\text{tree}]} \frac{(pq)^{-|E(\Gamma)|}}{|\aut (\Gamma)|} \sub(\Gamma,f^+),
\end{equation}
where $\Gamma$ is a tree in $S(\{x,y,z\})$ without edges connecting two vertices of type $x$ or two vertices of type $y$. 
\end{proposition}

\subsection{Tree Expansion of the Asymptotic Kontsevich Integral}

Let $\Gamma$ be a tree in $S(X)$ and for every vertex $v$ of $\Gamma$, denote by $y_v$ the cohomology class in $H^1(\Gamma_{\!\circ},\Z)$ of any edge of the circle inserted at $v$. 
Using Proposition \ref{subs}, one obtains the following formula:
$\sub(\Gamma,f^+)=\prod_v \frac{1}{(k_v-1)!}\frac{\partial^{k_v} f^+(v,y_v)}{\partial y_v^{k_v}}$.

Take formally $p=q$ in  formula \eqref{ntor} and apply the operator $p^{-2\deg}$. Denote by $F^{+(n)}$ the $n$-th derivative of $F^+$. Then, the following formula holds:

\begin{equation*}
\begin{split}
p^{-2\deg}Z^{\sun}_{>1}T(p,p)&=
\frac{\exp(\Theta/48)}{p^{-2\deg}\langle \Omega,\Omega\rangle}\sum_{[\Gamma,\text{tree}]}\frac{p^{-4|E(\Gamma)|}}{|\aut(\Gamma)|}\\
&\prod_{v\text{ type }z} \frac{(-p^2)^{k_v}F^{+(k_v)}(y_v)}{(k_v-1)!}
\prod_{v\text{ type }x,y} \frac{p^{k_v}F^{+(k_v)}(\frac{y_v}{p})}{(k_v-1)!}.
\end{split}
\end{equation*}
Using the equation $\sum_v k_v=2 |E(\Gamma)|$, the estimate $F^{+(n)}(h)=\frac{(-1)^{n-1}(n-1)!}{2h^n}+O(1)$ and letting $p$ go to infinity, one obtains our last proposition.

\begin{proposition}
The asymptotic Kontsevich integral may be computed by the following formula:
\begin{equation*}
\begin{split}
\lim_{p\to\infty}p^{-2\deg}Z^{\sun}_{>1}T(p,p)&=\\
\exp(\Theta/48)&\sum_{[\Gamma,\text{tree}]}\frac{1}{|\aut(\Gamma)|}
\prod_{v\text{ type }z} \frac{F^{+(k_v)}(y_v)}{(k_v-1)!}
\prod_{v\text{ type }x,y} \frac{-1}{2y_v^{k_v}},
\end{split}
\end{equation*}
where $\Gamma$ is a tree in $S(\{x,y,z\})$ without edges connecting two vertices of type $x$ or two vertices of type $y$.
\end{proposition}


\bigskip
\noindent
ETH Z\"urich, R\"amistrasse 101, CH-8092 Z\"urich, Switzerland.

\noindent
\texttt{sebastian.baader@math.ethz.ch}

\bigskip
\noindent
UPMC Universit\'e Paris 6, Institut de Math\'ematiques de Jussieu, F-75005, Paris, France.

\noindent
\texttt{marche@math.jussieu.fr}

\end{document}